\documentclass[twocolumn,aps,showpacs,amsmath,amssymb,floatfix]{revtex4}
\usepackage{graphicx}
\usepackage{array}
\usepackage{float}
\begin{document}
\title{Excited Random Walk in One Dimension}
\author{T. Antal}
\altaffiliation{On leave from Institute for Theoretical Physics -- HAS, E\"otv\"os University, Budapest, Hungary}
\affiliation{Center for Polymer Studies and Department of Physics,
Boston University, Boston, Massachusetts 02215}
\author{S.~Redner}\email{redner@cnls.lanl.gov}\altaffiliation{Permanent address: 
Department of Physics, Boston University, Boston, Massachusetts, 02215 USA}
\affiliation{Theory Division 
and Center for Nonlinear Studies, Los Alamos National Laboratory, Los
Alamos, New Mexico 87545 USA}

\begin{abstract}
  
  We study the excited random walk, in which a walk that is at a site that
  contains cookies eats one cookie and then hops to the right with
  probability $p$ and to the left with probability $q=1-p$.  If the walk hops
  onto an empty site, there is no bias.  For the 1-excited walk on the
  half-line (one cookie initially at each site), the probability of first
  returning to the starting point at time $t$ scales as $t^{-(2-p)}$.  Although
  the average return time to the origin is infinite for all $p$, the walk
  eats, on average, only a finite number of cookies until this first return
  when $p<1/2$.  For the infinite line, the probability distribution for the
  1-excited walk has an unusual anomaly at the origin.  The positions of the
  leftmost and rightmost uneaten cookies can be accurately estimated by
  probabilistic arguments and their corresponding distributions have
  power-law singularities near the origin.  The 2-excited walk on the
  infinite line exhibits peculiar features in the regime $p>3/4$, where the
  walk is transient, including a mean displacement that grows as $t^\nu$,
  with $\nu>\frac{1}{2}$ dependent on $p$, and a breakdown of scaling for the
  probability distribution of the walk.

\end{abstract}

\pacs{02.50.-r, 05.40.-a}

\maketitle

\section{Introduction}

In this work, we study the $k$-excited random walk in one dimension.  This
recently-introduced model \cite{PW97,D99}, gives rise to a stochastic process
with many unusual properties \cite{P01,ABV03,BW03,Z04}.  The definition of
the excited random walk is as follows: Initially each site of a
one-dimensional chain contains a fixed number of cookies $k$ (with $k\geq
1$).  A random walk makes nearest-neighbor hops between sites on the chain,
with hopping probabilities that depend on whether or not the current site
contains a cookie.  When the walker hops onto a site that contains one or
more cookies, the walker eats a cookie and, at the next time step, the walk
hops to the right with probability $p>1/2$ and to the left with probability
$q=1-p$.  If the walker hops onto an empty site, then the hopping
probabilities to the left and right at the next step are both 1/2.  This bias
can be viewed as the foraging strategy of a primitive organism.  When it
encounters food, the organism consumes the food and continues to move in a
preferred direction.  Conversely, the organism wanders aimlessly when it does
not encounter any food.

\begin{figure}[ht] 
 \vspace*{0.cm}
 \includegraphics*[width=0.13\textwidth]{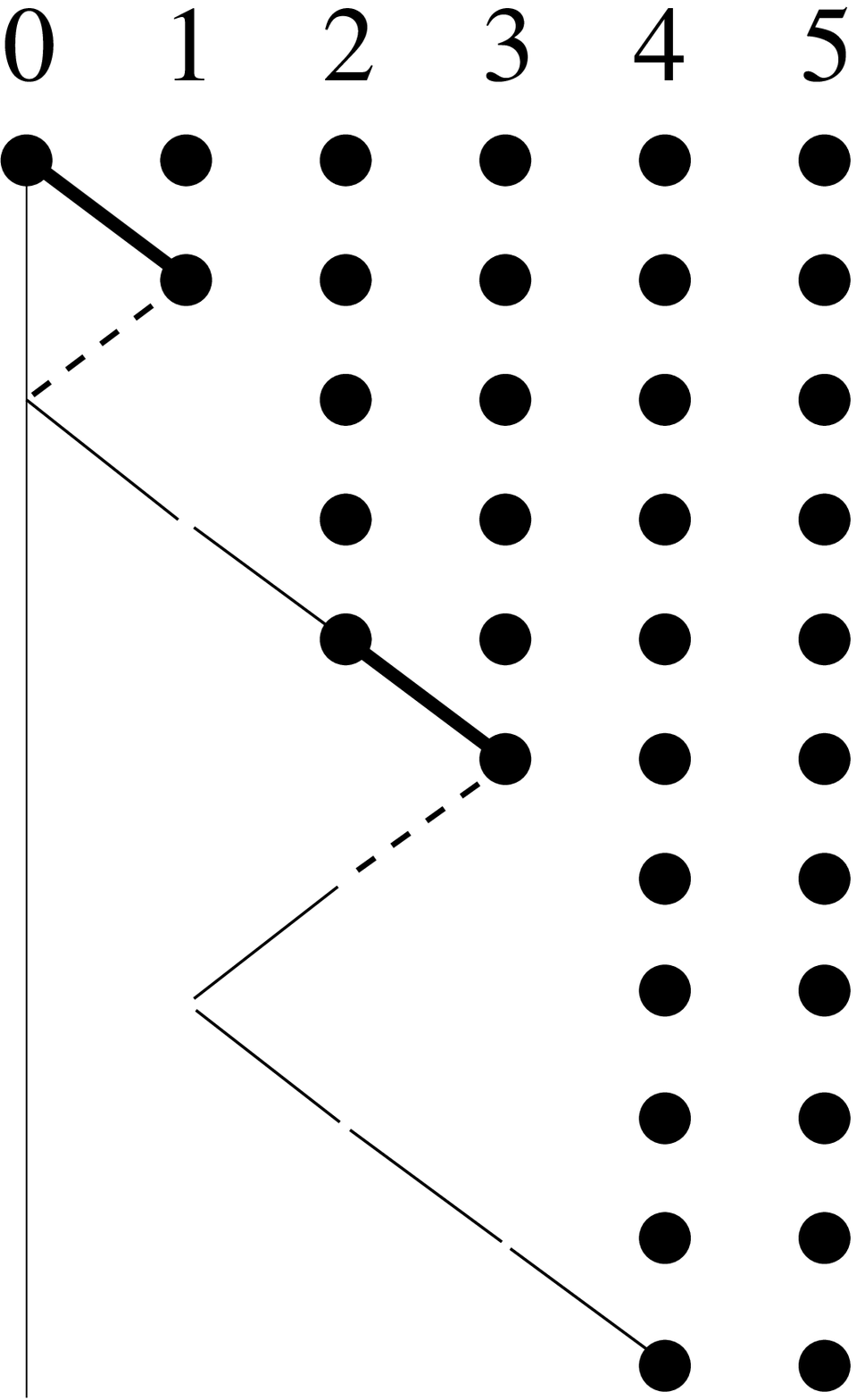} \hskip 0.5in \includegraphics*[width=0.25\textwidth]{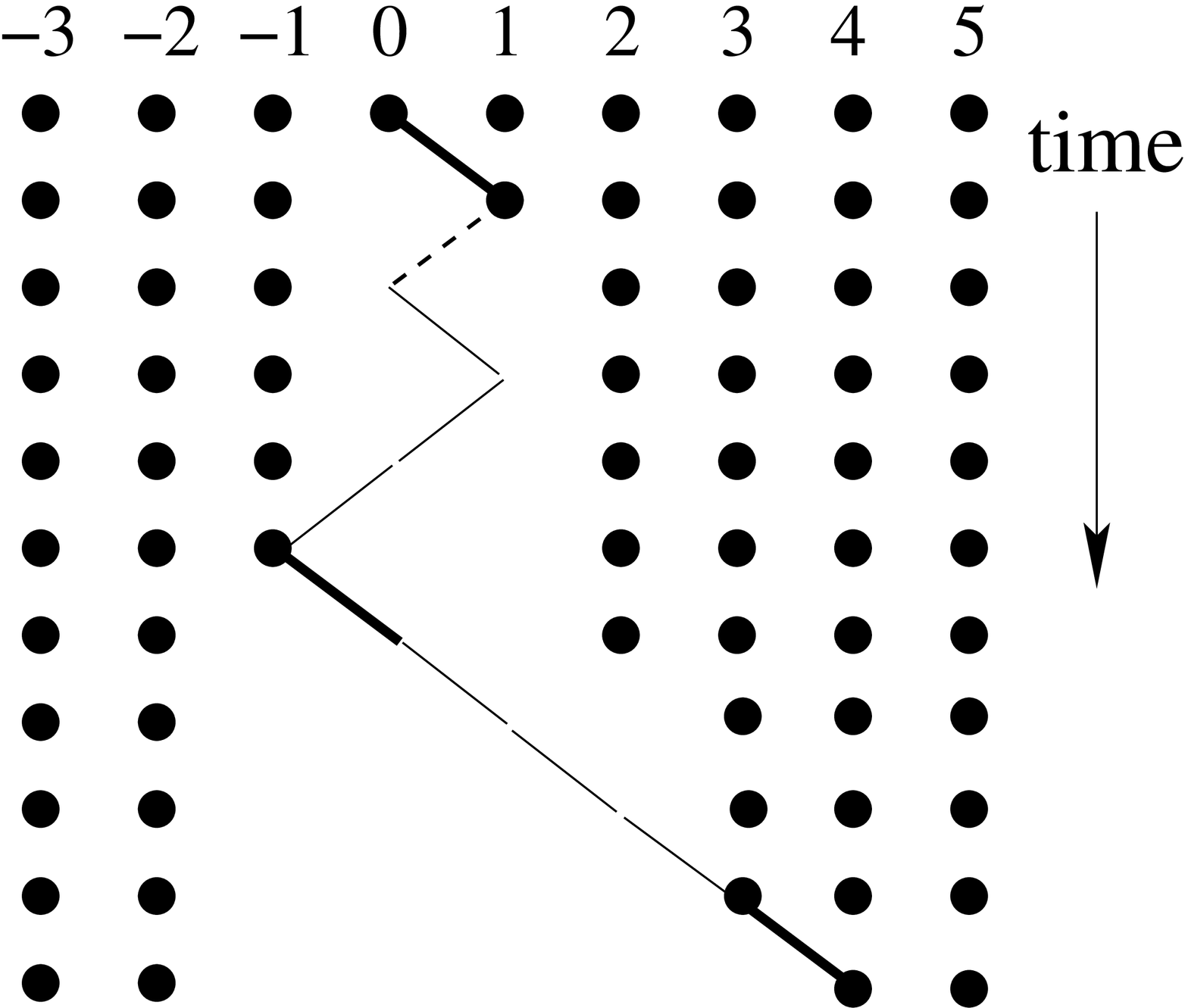}
 \caption{Space-time plot of the first few steps of a 1-excited random walk
   in one dimension with reflection at the origin (left), and on the infinite
   line (right).  Cookies are shown as large dots.  Steps occurring with
   probability $p$ and $q$ are shown as thick and dashed lines, while steps
   occurring with probability 1/2 are shown as thin lines.  For the half-line
   system, the last cookie at $n=4$ is first hit after 10 time steps.}
\label{def-1}
\end{figure}

A basic feature of the excited random walk is the transitory nature of the
bias; that is, when all the cookies at a given site are eaten the subsequent
motion is unbiased.  This non-Markovian feature is responsible for the
unusual properties of the excited random walk.  Related non-Markovian
constraints arise in many important stochastic phenomena.  Perhaps the best
known example is the self-avoiding random walk, where returns to
previously-visited sites are strictly forbidden \cite{SAW}.  This leads to an
expansion, relative to the pure random walk, in which the root-mean-square
displacement grows faster than $t^{1/2}$ when the spatial dimension is less
than 4.  A related example is the ``true'' self-avoiding walk, where returns
to previously-visited sites are inhibited by an energetic constraint
\cite{TSAW}.  The opposite type of perturbation occurs in the reinforced
random walk \cite{PV99}, where the walk has a preference to hop onto
previously-visited sites.  Depending of the functional form of this
reinforcement, the walk may, quite surprisingly, ``condense'' onto a support
that consists of only a finite number of sites as $t\to\infty$ \cite{PV99}.

The excited random walk represents a simple and surprisingly rich example of
a non-Markovian process in which the bias at each site disappears when all
the cookies on the site are eaten.  Because a random walk in one dimension is
recurrent, that is, the probability to eventually return to the starting
point equals one, the walk returns to previously-visited sites infinitely
often.  For the 1-excited walk, the reinforcement mechanism due to the
presence of cookies is insufficient to counterbalance this recurrence of the
random walk; thus the 1-excited random walk remains recurrent for any $p<1$.
Nevertheless, this model exhibits many novel features.  For example, the
probability that a 1-excited walk on an infinite line first returns to the
origin at time $t$ decays as $t^{-(2-p)}$, rather than the $t^{-3/2}$ decay
for a pure random walk \cite{DA01,DAA03,BW03}.

On the other hand, the $k$-excited random walk with $k>1$ is more strongly
influenced by the cookie-induced bias.  In fact, for $p>\frac{1}{2}
\left(1+\frac{1}{k}\right)$, the $k$-excited random walk becomes transient
\cite{Z04}.  This transition has dramatic consequences on the statistical
properties of the walk as $p$ is varied.  As a final introductory point, the
excited random walk appears to be most interesting in one dimension because
of the non-trivial competition between the reinforcement when a cookie is
eaten and the fact that encounters with sites that contain cookies become
progressively less frequent.  In contrast, an unbiased random walk in greater
than one dimension typically encounters previously-unvisited sites
\cite{F68,fpp}.  Because each such site contains a cookie, the reinforcement
mechanism acts as a non-zero global bias that is superimposed on the random
walk motion.  Thus the asymptotic behavior of an excited walk in higher
dimensions is close to that of a simple biased random walk.

In the next section, we begin by recounting some basic first-passage
properties for diffusion in a finite interval \cite{F68,fpp}.  We present
these results separately because they underlie our analytical treatment of
excited random walks.  In Sec.~III, we study a particularly simple version of
the excited random walk -- the 1-excited walk with a reflecting boundary
condition at the origin.  We use basic tools of first-passage processes to
determine the time dependence of the size of the region in which the cookies
have been eaten and thence the probability for the walk to first return to
the origin.  The long-time behavior of the first-passage probability in a
closely related model was obtained previously, using quite different methods,
in \cite{DA01,DAA03}.  Our approach is also complementary to the tools
employed in the probabilistic studies of excited random walks given in
Ref.~\cite{P01,ABV03,BW03,Z04}.  In Sec.~IV, we first present numerical
simulations for the unusual properties for the probability distribution and
for the span of the 1-excited random walk on the infinite line.  We also
construct a probabilistic argument to quantify the asymmetry in the number of
cookies eaten from the left and right sides of the line.  Finally, in Sec.~V
we present several intriguing properties of the 2-excited walk on the
infinite line, including the lack of scaling in the probability distribution,
the time for the 2-excited walk to traverse a finite domain, and basic
features of the cookie density profile that is left behind as the walk
advances.

\section{Preliminaries: First Passage in an Interval}

To determine the motion of excited random walks, we require the first-passage
properties of a nearest-neighbor unbiased random walk, or equivalently,
continuum diffusion, in a finite interval $[0,k]$ when both boundaries are
absorbing.  The first such property is the eventual hitting, or exit,
probability to $x=0$ and to $x=k$ for a random walk that starts at $x_0$
within the interval \cite{F68,fpp}.  These exit probabilities are:
\begin{equation} E_-(k,x_0)=1-\frac{x_0}{k}\qquad
  E_+(k,x_0)=\frac{x_0}{k}\,.
\label{exitprob_abs}
\end{equation}
Closely related quantities are the conditional exit times to $x=0$ and $x=k$.
These are the mean times for the random walk to hit $x=0$ or to hit $x=k$
with the restriction that the other boundary is never hit.  These conditional
times are:
\begin{equation}
t_-(k,x_0)=\frac{x_0}{3}(2k-x_0)\qquad t_+(k,x_0)=\frac{1}{3}(k^2-x_0^2) ~.
\label{exittime_abs}
\end{equation}
Next, we will also need the Laplace transforms of the first-passage
probabilities for a random walk that starts at $x_0$ to exit at the left edge
and at the right edge of the absorbing interval $[0,k]$ at time $t$.  In the
continuum limit, these first passage probabilities are (Eq.~(2.2.10) in
Ref.~\cite{fpp}):
\begin{eqnarray}
\label{lr-def}
\ell_k(x_0,s)&=&
\frac{\sinh\sqrt{2s}\,(k-x_0)}{\sinh\sqrt{2s}\,k}\,,\nonumber\\
r_k(x_0,s)&=&
\frac{\sinh\sqrt{2s}\,x_0}{\sinh\sqrt{2s}\,k}\,,
\end{eqnarray}
where we take the diffusion coefficient to be $D=1/2$ for the
nearest-neighbor random walk.

We will also need first-passage properties in the interval $[0,k]$ when the
boundary at $x=0$ is reflecting and the other boundary is absorbing.  We
implement the reflecting boundary condition for a discrete random walk that
is at site $x=0$ by allowing the walk to remain at $x=0$ with probability 1/2
or by having hop to $x=1$ with probability 1/2.  This definition for the
reflecting boundary has the advantage that it leads to a constant probability
density around this boundary for large times.  With this definition, the average
time to exit at the right edge of the interval $[0,k]$ as a function of $x_0$
is
\begin{equation}
t (k, x_0) = k(k+1)-x_0(x_0+1) ~.
\label{exittime_ref}
\end{equation}
Finally, the Laplace transform of the first-passage probability to exit at
the right edge at time $t$ in the continuum limit is (Eq.~(2.2.21) in
\cite{fpp})
\begin{eqnarray}
\label{R-def}
{\cal R}_k(x_0,s) =\frac{\cosh\sqrt{2s}\,(k-x_0)}{\cosh\sqrt{2s}\,k}.
\end{eqnarray}

We will exploit these results to now determine various basic features of the
1-excited random walk in one dimension.  There are two basic properties of
the excited random walk that we will be able to determine analytically: (i)
the location of the last (leftmost) cookie; we obtain this property in terms
of the first-passage statistics of a 1-excited random walk to the right edge
of the cookie-free domain, and (ii) the statistics of returns to the origin.

\section{1-Excited Random Walk on the Half-Line}

Consider the 1-excited random walk that is restricted to the non-negative
integers by imposing a reflecting boundary condition at the origin.
Initially there is one cookie at each site and the walk is at $n=0$.  At the
first time step, the walk either hops to $n=1$ with probability $p$ or it
remains at site $n=0$ with probability $q$ (Fig.~\ref{def-1}).  At later
times, a bias in the random walk motion occurs whenever the walk reaches the
right edge of the cookie-free domain.  The restriction of fixing one boundary
condition while the other boundary condition ``floats'' greatly simplifies
the analytical description of the 1-excited walk on the half-line compared to
the infinite line system.

\subsection{Time to Hit the Last Cookie}

It is useful to characterize the walk by the sequence of times $t_n$ at which
the walk hits the leftmost (last) cookie at $n$ (Fig.~\ref{def-1}).  Suppose
that the walk has just reached the cookie at $n-1$.  Then with probability
$p$, the next step of the walk is to right, which takes one time step, while
with probability $q$, the walk steps left to $n-2$ and the last cookie now is
at $n$.  In this latter case, the average time it takes the walk to reach
site $n$ is $t(n,n-2)=4n-2$ [see Eq.~(\ref{exittime_ref})] plus the time for
the initial step from $n-1$ to $n-2$.  Thus the mean time increment for the
last cookie to move from $n-1$ to $n$ is
\begin{eqnarray*}
\langle dt_n\rangle = p\cdot 1 +q\cdot (1+4n-2) = 1+q(4n-2) ~.
\end{eqnarray*}
Consequently, the mean time for the walk to hit the last cookie at $n$ for
the first time is
\begin{equation}
\label{tav}
\langle t_n\rangle =  \sum_{k=1}^n \langle dt_k\rangle = 2qn^2+n
\end{equation}
Notice that for $q\to 0$ ($p\to 1$), the early time behavior is governed by
ballistic motion, after which there is a crossover at a time of order $1/q$
to the asymptotic diffusive behavior of $\langle t_n\rangle \sim 2q n^2$.

We now extend this analysis to obtain the distribution of first-passage times
for the walk to hit the last cookie at $n$ at time $t$.  Let $Q_k(t)$ be the
probability that the time interval between the walk first hitting site $k-1$
and first hitting site $k$ equals $t$.  Similarly, let $F_n(t)$ be the
first-passage probability that the walk hits site $n$ for the first time at
time $t$.  Then
\begin{eqnarray*}
F_n(t)= \int \prod_{k=1}^n Q_k(t_k) dt_k\,,
\end{eqnarray*}
with $t_1+t_2+t_3+\cdots t_n=t$.  Because this relation is a convolution, the
Laplace transform $F_n(s)$ obeys
\begin{equation}
\label{Fns}
F_n(s) = \prod_{k=1}^{n} Q_k(s).
\end{equation}

Now the Laplace transform of $Q_k(t)$ is, using Eq.~(\ref{R-def}),
\begin{eqnarray}
Q_k(s) &=& p(1-s) + q(1-s){\cal R}_k(2,s)\,.
 \label{Qfull}
\end{eqnarray}
The first term accounts for the direct step from $k-1$ to $k$.  The
contribution of this direct process to the discrete generating function would
be $pz$.  In the continuum limit, the factor $z$ maps to the Laplace variable
$1-s$.  The second term accounts for a single step to the left (factor
$q(1-s)$) and the subsequent first-passage from $k-2$ to $k$ in the presence
of a reflecting boundary at the origin.  The first step of the walk is
accounted for by defining $Q_1(s)=e^s$.

We anticipate that the dominant contribution to $F_n(s)$ in the long-time
limit comes from terms where $k\sim\sqrt t$ \cite{AHZ}.  
In this limit,
Eq.~(\ref{Qfull}) reduces to
\begin{equation}
 Q_k(s) = 1 - 2q\sqrt{2s}\tanh\sqrt{2s}\,k ~,\nonumber
\end{equation}
by expanding ${\cal R}_k$ of Eq.~(\ref{R-def}) for large $k$.  Using this
approximation for $Q_k$ and taking the limit $s\to 0$, corresponding to
long-time behavior, $F_n(s)$ can be written as
\begin{eqnarray}
 F_n(s) &\approx& \prod_{k=1}^n \left(1-2q\sqrt{2s}\tanh\sqrt{2s}\,k\right)\,,\nonumber \\
&=& \exp\left(\sum_1^n \ln(1-2q\sqrt{2s}\tanh\sqrt{2s}\,k)\right)\,'\nonumber \\
&\sim& \exp\left(-\int^x 2q\sqrt{2s}\tanh\sqrt{2s}\,k\,\,dk\right)\,'\nonumber \\ 
&\sim& (\cosh\sqrt{2s}\,x)^{-2q} ~,
 \label{Fclosed}
\end{eqnarray}
where we use $\int \tanh x\,dx=\ln\cosh x$ and we replace the discrete
variable $n$ by the continuous variable $x$ in the last two lines.

{}From this formal solution, the cumulants $\kappa_\ell$ of the first-passage
time may be obtained as the coefficients of $(-s)^\ell$ in the power series
of $\ln F(x,s)$. For this we use the product representation of the $\cosh$
function, and the series expansion of $\ln(1+z)$ to give
\begin{eqnarray*}
\ln F(x,s) = 2q \sum_{\ell=1}^\infty \frac{1}{\ell} \left( 
\frac{8x^2}{\pi^2}\right)^\ell (-s)^\ell
\sum_{k=0}^\infty(2k+1)^{-2\ell} ~.
\end{eqnarray*}
Recognizing that the second sum can be expressed in terms of the Riemann zeta
function, the cumulants are
\begin{equation}
  \kappa_\ell =  qx^{2\ell} C(\ell) ~,
\end{equation}
where
\begin{eqnarray*}
  C(\ell) = 8^{\ell+1/3} (1-2^{-2\ell}) (\ell-1)!  \pi^{-2\ell}
\zeta(2\ell)\,,
\end{eqnarray*}
which gives $C(1)=2, C(2)=4/3, C(3)=32/15$, {\it etc}.  This result agrees
with the direct solution given in Eq.~(\ref{tav}) for the average
first-passage time $\langle t(x) \rangle = \kappa_1 = 2qx^2$.

Finally, the time dependence of the first-passage probability for the last
cookie -- the probability that the cookie at $x$ is first hit at time $t$ --
is obtained by inverting the Laplace transform in Eq.~(\ref{Fclosed}).  To
simplify this calculation, we introduce the scaling variables
\begin{eqnarray*}
 \mu = \frac{x}{\sqrt t} ~ ,~~~ y=st ~,
\end{eqnarray*}
to rewrite the first-passage probability $F(x,t)$ in the scaling form
\begin{equation}
 \label{Phidef}
  t F(x,t) \equiv \Phi(\mu) = \frac{1}{2\pi i} \int\limits_{-i\infty}^{i\infty} dy~ e^{y}  
  (\cosh\sqrt{2y}\mu)^{-2q} ~.
\end{equation}

\begin{figure}[htb]
  \centering
  \includegraphics[width=0.55\textwidth]{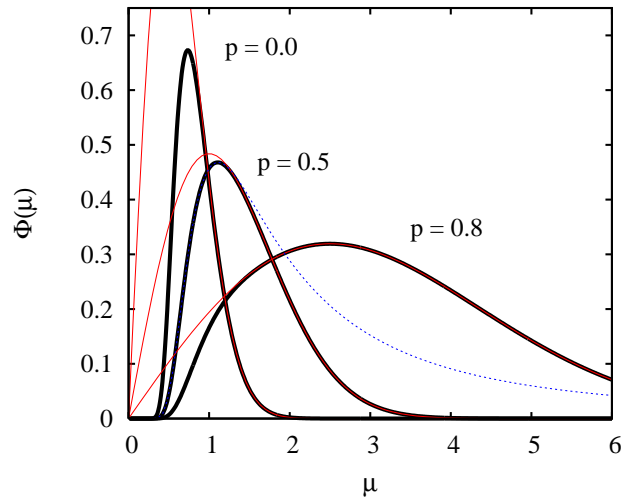} 
  \caption{The scaling function $\Phi(\mu)$ in the scaling limit for $p=0,
    0.5$, and $0.9$.  The large-$\mu$ asymptotic is represented by the first
    term in Eq.~(\ref{Philarge}) (thin curves).  For $p=0.5$ the first term
    of the small $\mu$ series of Eq.~(\ref{Phismall}) is also shown (dotted).
  }
  \label{fig:rphi}
\end{figure}

For the case of no bias ($p=q=1/2$), the integrand has simple poles on the
negative real axis at
\begin{eqnarray*}
 y = -\frac{8y\mu^2}{(2k+1)^2\pi^2} ~, ~~~ k=0,1,\ldots ~.
\end{eqnarray*}
By evaluating the corresponding residues at these poles, we find the scaling
function
\begin{equation}
 \label{Phismall}
 \Phi_{p=\frac{1}{2}}(\mu) = 
 \frac{\pi}{2\mu^2} \sum_{k=0}^\infty (2k+1)(-1)^k \,
e^{ -(2k+1)^2\pi^2/8\mu^2} ~.
\end{equation}
Because of the presence of an exponential cutoff, this series is easy to
compute numerically and the result is displayed in Fig.~\ref{fig:rphi},
together with the small-$\mu$ limiting behavior (corresponding to small $x$
for fixed $t$ or large $t$ for fixed $x$) that stems from the first term of
Eq.~(\ref{Phismall}):
\begin{eqnarray*}
 \Phi_{p=\frac{1}{2}}(\mu) \approx \frac{\pi}{2\mu^2} e^{-\frac{\pi^2}{8\mu^2}} ~.
\end{eqnarray*}

\begin{figure}[htb]
  \centering
  \includegraphics[width=0.55\textwidth]{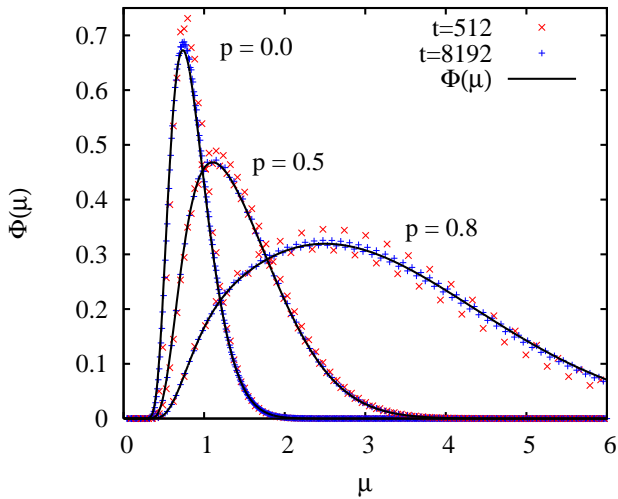} 
  \caption{The scaling function compared to the numerical exact propagation 
    of the probability distribution. }
  \label{fig:rnumer}
\end{figure}

For general $p\ne1/2$ the poles in Eq.~(\ref{Phidef}) are branch points and
we therefore evaluate this integral in a simpler manner by first using the
binomial theorem to write Eq.~(\ref{Fclosed}) as
\begin{equation}
 F(x,s) = 2^{2q} \sum_{k=0}^\infty \binom{-2q}{k}  e^{-2x\sqrt{2s}\,(k+q)} ~.
\end{equation}
Each term in this series can be easily back transformed term by term (see
29.3.82 in Ref.~\cite{AS}).  The scaling function thus has the form
\begin{equation}
 \Phi(\mu) = \mu \frac{2^{2q+1/2}}{\sqrt{\pi}} \sum_{k=0}^\infty 
\binom{-2q}{k} (k+q) \, e^{-2\mu^2 (k+q)^2} ~.
 \label{Philarge}
\end{equation}
For $q=1/2$, this function is identical to $\Phi(\mu)$ in
Eq.~(\ref{Phismall}) in spite of the apparent visual differences in these two
formulae, since both functions originate from the same Laplace transform.
Eq.~(\ref{Philarge}) can be viewed as the large-$\mu$ expansion of
$\Phi(\mu)$.  As shown in Fig.~\ref{fig:rphi}, the $k=0$ term in this series
representation
\begin{equation}
 \Phi(\mu) \approx \mu q \frac{2^{2q+1/2}}{\sqrt{\pi}} e^{-2(\mu q)^2},
\end{equation}
together with the $k=0$ term in the small-$\mu$ series for
Eq.~(\ref{Phismall}) give a fairly good approximation for $\Phi(\mu)$ over
the entire range of $\mu$ for $p=1/2$.  Fig.~\ref{fig:rnumer} shows the
solution Eq.~(\ref{Philarge}), together with the exact numerical propagation
of the probability distribution at different times (see Appendix).

\subsection{Location of the Last Cookie}

Another important characterization of the excited random walk is the
probability distribution for the location of the last cookie.  Now the last
cookie either remains stationary or it moves only to the right during
each step of the underlying excited random walk.  Consequently, we can write
a simple relation between the first-passage probability for hitting the last
cookie, derived in the previous section, and the probability distribution of
this last cookie.

The probability that the last cookie is located at site $n$ at time $t$
equals the joint probability that a walk first hopped to site $n-1$ {\it
  before} time $t$, and that this walk then first hopped to site $n$ {\it
  after} time $t$.  That is,
\begin{equation}
 P_n(t) = \sum_{t'=1}^t \sum_{t''=t-t'+1}^\infty F_{n-1}(t')\, Q_{n}(t'')  ~.
\end{equation}
In the continuum limit, the Laplace transform of this equation is
\begin{equation*}
P_n(s) = \frac{1}{s}F_{n-1}(s)\left[ 1-Q_n(s) \right] = 
\frac{1}{s}\left[ F_{n-1}(s) - F_n(s)\right] ~,
\end{equation*}
where we used Eq.~(\ref{Fns}) to simplify the last factor.  In the continuous
space limit, the relation between $P$ and $F$ becomes
\begin{equation}
 P(x, s) = - \frac{1}{s} \frac{\partial F(x,s)}{\partial x} ~,\nonumber
\end{equation}
which is equivalent to
\begin{equation}
\label{P-F}
 P(x,t) = -\frac{\partial}{\partial x} \int_0^t dt' F(x,t') ~.
\end{equation}
Substituting the scaling form for $F(x,t)$ into Eq.~(\ref{P-F}) and
performing the derivative, we express the integrand in scaled units, after
which the integral can be evaluated easily.

These steps lead to the following simple connection between the first-passage
and occupation probabilities of the last cookie:
\begin{equation}
F(x,t)=\frac{1}{t} \Phi(\frac{x}{\sqrt t})\,,\quad{\rm and}\quad 
P(x,t)=\frac{2}{x} \Phi(\frac{x}{\sqrt t}) ~,
\label{scaling}
\end{equation}
with the {\it same} scaling function for both quantities.  We exploit this
connection between $F$ and $P$ to obtain the finite time data in
Fig.~\ref{fig:rnumer} by first calculating the position distribution of the
last cookie from the numerical evaluation of the Master equation (see
Appendix), and then using Eq.~(\ref{scaling}) to arrive to the scaling
function $\Phi(\mu)$.

Finally, we use Eqs.~(\ref{Philarge}) and (\ref{scaling}) to obtain the exact
expression for the average position of the last cookie:
\begin{equation}
  \label{xaver}
  \langle x(t) \rangle = \int_0^\infty dx~ x P(x,t) \equiv 2\sqrt{t}\,{\cal N}(q) ~,
\end{equation}
where the proportionality constant
\begin{equation}
{\cal N}(q) = \int_0^\infty d\mu~ \Phi(\mu) = \frac{2^{2q-3/2}}{\sqrt{\pi}\,q}
  \frac{\left[\Gamma(1+q)\right]^2}{\Gamma(1+2q)} \nonumber
  \label{Phinorm}
\end{equation}
is obtained by first integrating (\ref{Philarge}) term by term and then using
the identity \cite{knuth}
\begin{eqnarray*}
\sum_{k=0}^\infty \binom{-2q}{k} \frac{1}{k+q} = 
\frac{\left[\Gamma(1+q)\right]^2}{q\Gamma(1+2q)} ~.
\end{eqnarray*}
The surprising point of this calculation is that $\langle x(t)\rangle\sim
\sqrt{t}/q$ as $q\to 0$, while the mean first-passage time for the walk to
reach a distance $x$ scales as $\langle t(x)\rangle\propto qx^2$.  More
precisely,
\begin{equation}
  \langle t(x) \rangle = 2q x^2\,,\qquad \langle x(t) \rangle \sim
 \frac{1}{\sqrt{2\pi}\,q} \sqrt{t} ~,
 \label{weird}
\end{equation}
while naive dimensional scaling would suggest that $\langle t\rangle\propto
qx^2$ implies $\langle x(t)\rangle\propto \sqrt{t/q}$.  The fact that this
naive scaling does not hold is a result of the broadness of the distribution
of first-passage times for the last cookie.

\subsection{First Passage to the Origin}
\label{origin}

As a complement to the statistical properties of the last cookie, we now
study first passage to the origin, namely, the probability for the 1-excited
walk to reach $x=0$ for the first time at time $t$.  We first provide a
heuristic argument (previously given in Ref.~\cite{BW03}) that the 1-excited
walk is recurrent and then we extend this line of reasoning to argue that the
first-passage probability asymptotically scales at $t^{-q-1}$.

If the walk is at a point $x-1>0$ for the first time, then the probability
for the walk to subsequently reach $x$ without hitting the origin is
$p+q\left(1-\frac{2}{x}\right) = 1-\frac{2q}{x}$.  Here the factor $p$
accounts for the event where the walk immediately eats a cookie by making a
step to the right, and the factor $q\left(1-\frac{2}{x}\right)$ [see
Eq.~(\ref{exitprob_abs})] accounts for the probability that the first step is
to the left and that the walk eventually eats another cookie on the right
without hitting the origin. Thus, if the walk starts at site $x=1$ with the
system initially full of cookies, the probability for the walk to reach (at
least) an arbitrary point $x>1$ without hitting the origin is
\begin{eqnarray}
\label{P-prod}
P(x)&=&\prod_{k=2}^x \left(1-\frac{2q}{k}\right) 
= \exp\left[\sum_{k=2}^x \ln\left(1-\frac{2q}{k}\right)\right]\nonumber\\
&\sim& \exp\left(-\int^x \frac{2q}{y}\, dy\right)\nonumber \\
&\sim& x^{-2q}.
\end{eqnarray}
Since $x^{-2q}\to 0$ as $x\to\infty$, the random walk must necessarily return
to the origin; that is, the 1-excited random walk is recurrent, unless $p=1$.

Since $P(x)$ is the probability that the walk eats {\it at least} $x$ cookies
without hitting the origin, $P'(x)\sim x^{-1-2q}$ is the probability that the
walk eats exactly $k$ cookies.  The average number cookies eaten by a walk
that does not hit the origin is therefore $\langle n \rangle\sim\int^\infty
dx~x^{-2q}$.  Thus for $q>1/2$ $(p<1/2$) a walk that manages to avoid the
origin typically wanders only a finite distance away (for a more rigorous
argument see Sec.~\ref{asym}).

If the walk is at site $x-1$ for the first time, the conditional average time
to reach site $x$ without hitting the origin is one with probability $p$, and
$\frac{4}{3}(x-1)$ with probability $q$ (see Eq.~(\ref{exittime_abs})), that
is $\langle dt_x\rangle=1+\frac{4q}{3}(x-1)$.  Thus the average time to reach
site $x$, without hitting the origin, is $\langle t_x \rangle =
1+\sum_{k=3}^x \langle dt_k\rangle\sim x^2$ for large $x$.  From this
relation between $t$ and $x$, we infer that the survival probability, namely,
the probability that the walk does not reach the origin by time $t$, should
vary as $S(t)=P(x)\sim x^{-2q}\sim t^{-q}$.  Consequently, the first-passage
probability $F(t)=-\frac{\partial S}{\partial t}$ should asymptotically decay
as
\begin{equation}
\label{F-state}
F(t)\sim t^{-q-1}\,.
\end{equation}
For $q<1$, the average time to reach the origin is infinite.  Quite
strangely, although it takes infinite time for the walk to hit the origin for
any $q$, the walk eats only a finite number of cookies, {\it i.e.}, wanders
only a finite distance from the origin, for $q>1/2$.  This result remains
valid even if there is an initial cookie-free region at the origin, as this
initial condition would affect the pre-asymptotic factors in
Eq.~(\ref{P-prod}).

To compute the first-passage probability to the origin more precisely, it is
useful to define two auxiliary quantities, $L_k(t)$ and $R_k(t)$, as follows:
\begin{itemize}
\item $L_k(t)\!\equiv\!$ the probability of first hitting $k=0$ without hitting
  $k+1$ a time $t$ after the walk first hits $k$.
\item $R_k(t)\!\equiv\!$ the probability of first hitting $k+1$ without
  hitting $k=0$ a time $t$ after the walk first hits $k$.
\end{itemize}

In the Laplace domain, these first-passage probabilities are, using
Eq.~(\ref{lr-def}),
\begin{eqnarray*}
L_k&=& q(1-s) \ell_{k+1}(k-1,s)\,, \\
R_k&=& p(1-s) + q(1-s) r_{k+1}(k-1,s)\,.
\end{eqnarray*}
For example, for $L_k$, the first step must be to the left (factor $q(1-s)$)
and then the factor $\ell_{k+1}$ is the first-passage probability to the
origin, without hitting $k+1$, for a walk that starts at $k-1$.  Similarly,
for $R_k$, the term $p(1-s)$ accounts for the transition $k\to k+1$ in a
single step.  The factor $q(1-s)$ in the second term accounts for a single
step to the left after which the walk is at $k-1$ while the last cookie is at
$k+1$.  Then the factor $r_{k+1}(k-1,s)$ gives the first-passage probability
to $k+1$ without hitting the origin, when starting from $k-1$.

Finally, we obtain the Laplace transform of the first-passage probability to
the origin for a walk that begins at $x=1$ with the system initially full of
cookies by summing over all paths that contain $0, 1, 2, \ldots$ first
passages to the last cookie before the origin is reached.  This gives
\begin{eqnarray}
\label{Fs}
F(s)&\!=\!&L_1+R_1\Big(L_2+R_2\big(L_3+R_3(L_4+\ldots\nonumber \\
&\!=\!& L_1+R_1L_2+R_1R_2L_3+R_1R_2R_3L_4+\ldots
\end{eqnarray}

In the limit $s=0$, this Laplace transform is just the integral of the
first-passage probability over all time or, equivalently, the probability of
eventual return to the origin.  It is easy to verify that this return
probability equals one.  In the $s\to 0$ limit, the auxiliary probabilities
$R_k$ and $L_k$ are
\begin{equation}
\label{RLlimit}
R_k=p+q\frac{k-1}{k+1}=1-\frac{2q}{k+1};\quad L_k=\frac{2q}{k+1}=1-R_k.
\end{equation}
Then
\begin{eqnarray}
F(s=0)&=&1-R_1+R_1(1-R_2)+R_1R_2(1-R_3)\nonumber \\
      &{ }&~~~~~~~~~ +R_1R_2R_3(1-R_4)+\ldots~.
\end{eqnarray}
This expression obviously equals one for all $p<1$, while for extreme case of
$p=1$, $F(s=0)=0$.  Thus eventual return to the origin is certain for any
bias $p<1$.

More formally, we denote the sum of the first $m$ terms of Eq.~(\ref{Fs}) as
\begin{equation}
\label{Fms}
F_m(s) = \sum_{k=1}^m L_k \prod_{j=1}^{k-1}R_j ~.
\end{equation}
For finite $m$, we use Eq.~(\ref{RLlimit}) to give, in the limit $s\to 0$,
\begin{equation}
\label{F0s}
F_m(s=0) = \frac{2q}{\Gamma(2-2q)} \sum_{k=1}^m \frac{\Gamma(k+1-2q)}{\Gamma(k+2)} ~, 
\end{equation}
which goes to one as $m\to\infty$ \cite{concmat}.  The use of the $s=0$ forms
for $R_k$ and $L_k$ from (\ref{RLlimit}) remains valid for $m$ of the order
of $1/\sqrt{s}$ or less.  We now assume that the main contribution to $F(s)$
comes from terms with $m<1/\sqrt{s}$ in Eq.~(\ref{Fms}).  That is,
\begin{equation}
F(s) \approx \frac{2q}{\Gamma(2-2q)} \sum_{k=1}^{1/\sqrt{s}} 
\frac{\Gamma(k+1-2q)}{\Gamma(k+2)} ~.
\end{equation}
This finite sum can be written as the difference of the infinite sum (which
equals 1) and the sum from $1/\sqrt{s}$ to $\infty$.  In this second sum,
$k\gg1$ for all terms as $s\to\infty$.  Therefore we can replace the ratio of
gamma functions by $k^{a-b}$ and the sum by an integral.  These steps lead to
\begin{equation}
1-F(s) \sim \int_{1/\sqrt{s}}^\infty dk~ k^{-2q-1} \sim s^q ~
\end{equation}
This result for $F(s)$ implies that for large times the first-passage
probability to the origin indeed behaves as $F(t)\sim t^{-1-q}$.

\section{1-Excited Walk on the Infinite Line}

\subsection{First-Passage Probabilities}

Using the same first-passage analysis as that given for the 1-excited walk on
the half-line, we can immediately deduce the first-passage probabilities for
the 1-excited walk on the infinite line to hit the left side and the right
side of a cookie-free region.  Thus the first-passage probability for the
walk to hit the left side of this gap is the same as the first-passage
probability to the origin for the half-line system.  From Eq~(\ref{F-state}),
$F_L(t)\sim t^{-1-q}$, where we now write the subscript $L$ to denote hitting
the left side of the gap.  The first-passage probability to the right side of
the gap can immediately be obtained merely by interchanging $p$ and $q$; thus
$F_R(t)\sim t^{-1-p}$.  Since the exponents in both of these first-passage
probabilities are smaller than 2, the mean time to hit either edge of the
cookie-free domain is infinite.

\subsection{The Probability Distribution}

Another basic characteristic of the 1-excited walk is the probability
distribution of displacements.  An analytical treatment of this quantity is
more involved than the semi-infinite system because of the existence of two
moving boundaries.  We thus turn to numerical simulations to provide a
qualitative picture of this distribution.

\begin{figure}[ht] 
 \vspace*{0.cm}
 \includegraphics*[width=0.45\textwidth]{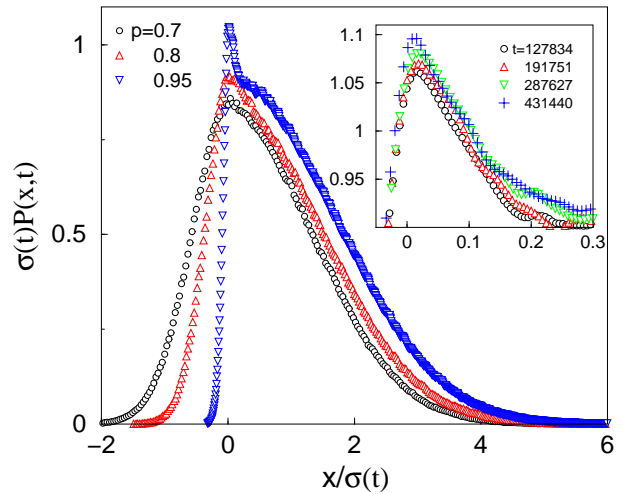}
 \caption{Scaled probability distributions of the 1-excited random walk for
   $2\times 10^6$ realizations for $ 85,223\approx 1.5^{28}$ time steps for
   $p=0.7, 0.8$ and 0.95.  Plotted on the horizontal axis is $x/\sigma$ and
   on the vertical axis is $\sigma\, P(x,t)$, where $P(x,t)$ is the
   probability distribution of the walk.  The data have been averaged over a
   $\frac{1}{2}$\% range.  The inset shows the distribution for $p=0.95$ near
   the origin for still longer times: $t=1.5^k$, with $k=29 - 32$ (127,834 --
   431,440 steps).  }
\label{pdf-1e}
\end{figure}

Fig.~\ref{pdf-1e} shows the scaled probability distribution for $p=0.7$, 0.8,
and 0.9.  For $p\agt 0.5$, the probability distribution has a nearly
Gaussian, but slightly asymmetric, appearance.  As $p$ increases toward 1,
this asymmetry becomes more pronounced, as shown in the figure.  In addition,
an apparent singularity in the distribution develops near the origin as $p\to
1$.  By running simulations for $p=0.7$ and $0.95$ to $431,440\approx
1.5^{32}$ steps, we observe quite accurate data collapse for $p=0.7$.
However for $p=0.95$, there is a small but persistent departure from data
collapse for $0\alt x/\sigma\alt 0.3$ (inset to Fig.~\ref{pdf-1e}.  We do not
understand whether this anomaly arises from a slow correction to scaling term
in the probability distribution or whether there exists an effect that lies
outside of scaling.

For all values of $p>1/2$, both $\langle x(t)\rangle$ and the variance
$\sigma(t)=\sqrt{\langle x(t)^2\rangle - \langle x(t)\rangle^2}$ grow as
$t^{1/2}$.  The time dependence for the displacement is easy to understand on
heuristic grounds.  Once a cookie-free domain of linear dimension $L$ has
been created, the domain then grows only when the walk reaches either edge of
the domain.  Since the motion of the walk is purely random within the domain,
the probability that the walk hits either edge is proportional to $1/L$.  If
the walk hits the right edge of the domain, then the walk will move
ballistically for typically $1/(1-p)=1/q$ steps as $p\to 1$.  Thus the
equation of motion for the domain length will be $\dot L= 1/(Lq)$ as $p\to
1$.  The solution to this equation gives $L\propto t^{1/2}$ for any value of
$p<1$.  However, it bears mention that this naive approach gives the wrong
value for the amplitude -- $1/\sqrt{q}$ -- rather than $1/q$ that was derived
in Eq.~(\ref{weird}).

\subsection{Asymmetry}
\label{asym}

A useful geometric characteristic of the 1-excited walk is the asymmetry in
the span of the walk.  Here the span is the length of the cookie-free region.
For an unbiased random walk, the growth of the span with time and its
distribution are well understood \cite{W}, and we now investigate how the
bias inherent in the 1-excited walk influences the span.  It is useful to
characterize the span by the positions of the last cookie on the right and
the last cookie on the left of the cookie-free region.  We can understand the
motion of these extreme cookies by the following heuristic argument.
Consider the situation where the cookie-free gap has length $L-2$ and where
the walk has just stepped to the extreme cookie on the right.  As shown at
the beginning of Sec.~\ref{origin}, the probability that the walk eats
another right cookie without eating a left cookie is $1-\frac{2q}{L}$.

With this basic result, the probability that the walk eats {\it precisely} $r>0$ consecutive
cookies (we term this event a single ``meal'') from the right edge of the
cookie-free region is
\begin{equation}
\begin{split}
P(r) &= \prod_{j=0}^{r-2}\left(1-\frac{2q}{L+j}\right) \times \frac{2q}{L+r-1}\\
&= 2q \frac{\Gamma(L)}{\Gamma(L-2q)} \frac{\Gamma(L+r-1-2q)}{\Gamma(L+r)} ~,
\end{split}
\label{prdef}
\end{equation}
where we use the convention $\prod_a^b=1$ for $a>b$. The probability that the walk eats precisely one cookie (the one at the starting position) is then $2q/L$.

When the initial cookie-free gap is large, this probability simplifies to
\begin{eqnarray*}
P(r) \sim 2q\,L^{2q}\,(L+r)^{-2q-1} ~.
\end{eqnarray*}
The $L$ dependence can be scaled out using the dimensionless variable $\tilde
r=r/L$ and the rescaled probability
\begin{equation}
\tilde P(\tilde r) \equiv LP(r) = 2q (1+\tilde r)^{-2q-1} ~.
\end{equation}
Note that this probability is normalized, $\int_0^\infty d\tilde r \tilde
P(\tilde r)$=1.  Using this probability, the average relative number of
consecutive cookies eaten from the right side of the gap, $\int_0^\infty
d\tilde r~ \tilde r \tilde P(\tilde r)$, is
\begin{equation}
\langle \tilde r \rangle= \langle r/L \rangle = \left\{ 
 \begin{split}
 \frac{1}{2q-1},  ~~~ q>1/2; \\
 \infty ,  ~~~ q\le1/2.
 \end{split} \right. 
\label{tilder}
\end{equation}
This result is paradoxical at first sight because, as we showed in
Sec.~\ref{origin}, the average time to return to the left side of the gap is
infinite.  However, the walk eats, on average, only a finite number of
cookies if $q\le1/2$, even though the time taken for its meal on the right
side of the gap is infinite.

The average relative number of eaten cookies from the left $\langle \tilde l
\rangle$ is also given by Eq.~(\ref{tilder}) after interchanging $q$ and $p$.
Although the 1-excited walk eats, on average, infinitely many consecutive
cookies from a single side of the gap during a single meal, the average ratio
$\langle r/(l+r) \rangle$ turns out to be finite for any value of $q$.  To
show this we compute the joint probability of eating $r$ consecutive cookies
on the right and then $l$ consecutive cookies on the left before the next
right cookie is eaten, under the initial condition that the walk starts at the
position of the right extreme cookie.  This joint probability can be derived
analogously to Eq.~(\ref{prdef}) 
\begin{equation}
\begin{split}
P(r,l) &=4pq\frac{\Gamma(L)}{\Gamma(L-2p)}\times\\
&\frac{\Gamma(L+l-1-2p)}{\Gamma(L+l-2q)}
\frac{\Gamma(L+l+r-1-2q)}{\Gamma(L+l+r)}\\
&\sim 4q\,(1-q)L^{2q}(L+r)^{1-4q} (L+r+l)^{2q-3}.
\end{split}
\end{equation}
It is again useful to define the rescaled variables $\tilde l=l/L$ and
$\tilde r=r/L$, as well as the rescaled probability
\begin{equation}
\tilde P(\tilde r, \tilde s) \equiv \tau^2 P(r,l) = 
4q(1-q) (1+\tilde r)^{1-4q} (1+\tilde s)^{2q-3} ~,
\label{lrdist}
\end{equation}
where $\tilde s=\tilde l+\tilde r$ is the total number of cookies eaten in
one right-left cycle.  One may check that this probability is normalized,
$\int_0^\infty ds \int_0^s dr \tilde P(r, s) = 1$.  

\begin{figure}[htb]
  \centering
  \includegraphics[width=0.5\textwidth]{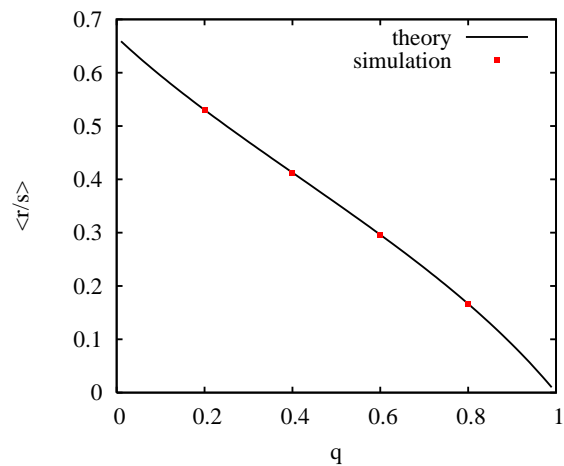} 
  \caption{Ratio of the number of right cookies eaten to total cookies eaten 
    after a single cycle (\ref{oneloop}). The squares are obtained from
    simulations with an initial gap ranging from 20 to 200.}
  \label{fig:lscomp}
\end{figure}

Finally, the average ratio for the number of right cookies to total cookies
eaten in a single cycle is
\begin{eqnarray}
\label{lovers}
\left\langle \frac{r}{s} \right\rangle &\!=\!& \left\langle \frac{\tilde r}{\tilde s} \right\rangle 
= \int_0^\infty ds \int_0^s dr ~\frac{r}{s}~\tilde P(r, s)\nonumber \\
&\!=\! &4q(1\!-\!q)  \int_0^\infty \frac{ds}{s}\,(1\!+\!s)^{2q\!-\!3} 
\int_0^s dr\, r(1\!+\!r)^{1\!-\!4q}. \nonumber \\
\end{eqnarray}
This integral can be performed with the help of 3.235 in Ref.~\cite{GR} (both
for $q<3/4$ and $q>3/4$), and also using some basic properties of the Euler
$\psi$ function \cite{AS} to give
\begin{equation}
\left\langle \frac{r}{s} \right\rangle
= 1 - \frac{2q^2}{(1-2q)^2} - \frac{2\pi q (1-q)}{(1-2q)(3-4q)} \cot 2\pi q
\label{oneloop}
\end{equation}
This gives the asymptotically exact ratio for the number of right cookies to
total cookies eaten in one cycle, under the assumption that the initial
cookie-free gap is large.  Although $\langle l \rangle$ or $\langle r
\rangle$ is infinite, the $q$-dependence of the ratio $\langle
r/(l+r)\rangle$ is well-behaved (Fig.~\ref{fig:lscomp}).

Despite the fact that the average number of consecutive cookies eaten from
the left or the right side of the gap in one cycle is infinite, our data
clearly show that the positions of the two extreme cookies both grow as
$t^{1/2}$, but with very different amplitudes. To obtain an approximation for
these amplitudes we now calculate the {\it typical} number $l^*$ and $r^*$ of
eaten consecutive cookies from each side.

\begin{figure}[ht] 
 \vspace*{0.cm}
 \includegraphics*[width=0.4\textwidth]{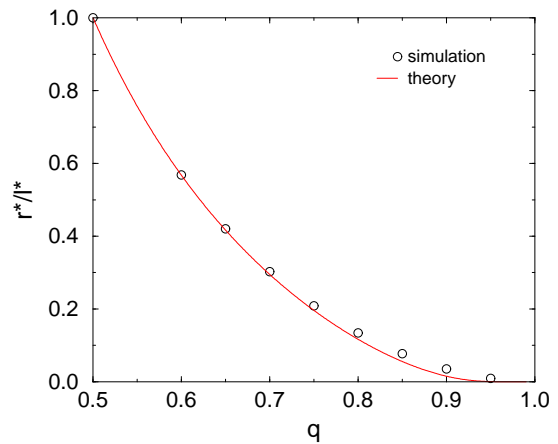}
 \caption{The ratio $r^*/l^*$ as a function of $p$ for the 1-excited random
   walk.  The smooth curve is the the prediction from Eq.~(\ref{lr-ratio}),
   while the data are based on $2\times 10^6$ realizations of the walk.}
\label{lr-compare}
\end{figure}

We estimate the value of $r^*$ by setting the probability of eating {\it at
  least} $r^*$ consecutive cookies on the right equal the probability of
going to the left before $r^*$ right cookies are eaten.  That is,
\begin{equation}
 \prod_{j=0}^{r^*-2}\left(1-\frac{2q}{L+j}\right) = \frac{1}{2} ~.
 \label{rstar}
\end{equation}
This criterion gives a crude estimate for the number of right cookies eaten
during a single consecutive string.  For $L\gg1$, Eq.~(\ref{rstar}) becomes
\begin{eqnarray*}
\left(\frac{L}{L+r^*}\right)^{2q} = \frac{1}{2} ~,
\end{eqnarray*}
which leads to the estimate that the typical length of the string of eaten right
cookies is $r^*=(2^{1/2q}-1)L$.  Applying this same argument to the left side
of the cookie-free region, a string of $l^*=(2^{1/2p}-1)L$ consecutive left
cookies are typically eaten. Their ratio thus equals the constant value
\begin{equation}
\label{lr-ratio}
\frac{r^*}{l^*}=\frac{2^{1/2q}-1}{2^{1/2p}-1}.
\end{equation}
As a function of time, $r^*/l^*$ clearly approaches a constant value as
$t\to\infty$ and it is easy to estimate its asymptotic value which is shown
in Fig.~\ref{lr-compare}, together with the prediction of
Eq.~(\ref{lr-ratio}).

\begin{figure}[ht] 
 \vspace*{0.cm}
 \includegraphics*[width=0.42\textwidth]{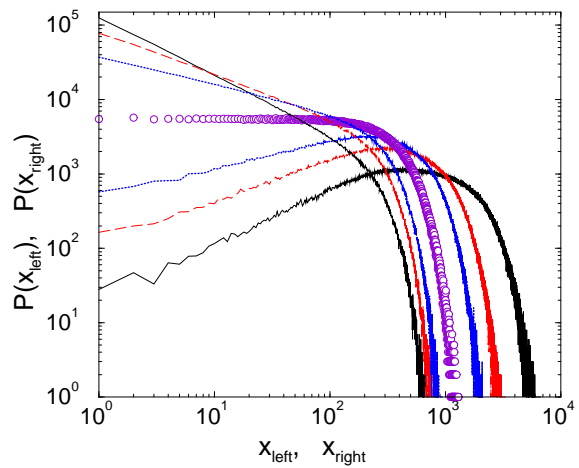}
 \caption{Probability distributions for the position of the rightmost cookie
   (positive slopes for small $x$) and the leftmost cookie for the 1-excited
   random walk.  Data are based on $2\times 10^6$ realizations for $85,223$
   time steps, for $p=0.7$ (line), $0.8$ (dashed), and $0.9$ (dotted).  For
   the rightmost cookie, the apparent power law for small $x$ has exponent
   0.35, 0.54, 0.74 and 0.84 for p=0.7, 0.8, 0.9 and 0.95, respectively,
   while for the leftmost cookie, the respective exponents are $-0.42$,
   $-0.63$, $-0.83$, and $-0.94$.  Also shown is the data for the unbiased
   random walk ($\circ$).}
\label{lr-cookie}
\end{figure}

To conclude this section, we point out another intriguing feature of the
extreme cookies; namely the probability distribution in their positions.
Quite strikingly, these distributions have short-distance power-law tails,
both for the left and for the right extreme cookies (Fig.~\ref{lr-cookie})
for $p\ne \frac{1}{2}$.  In contrast, for the pure random walk case of
$p=\frac{1}{2}$, the extreme cookie distribution is essentially the same as
the span distribution, which is known to be a Gaussian-like function
\cite{W}.  The origin of the short-distance power-law tails in the
probability distributions of the last cookie seem to be connected with the
anomalous scaling of the first-passage probability, but we do not have a
definitive explanation for this phenomenon.

\section{2-Excited Random Walk}

In the 2-excited random walk, each lattice site initially has two cookies and
the $(p,q)$ bias occurs as long as the walk lands on a site that contains at
least one cookie.  This 2-excited random walk undergoes a transition between
recurrence and transience as $p$ is increases beyond a critical value
$p_c=\frac{3}{4}$ \cite{Z04}. Indeed, the $k$-excited walk for any $k>1$
appears to have qualitatively similar behavior as the 2-excited walk.  For
simplicity, we thus focus on the 2-excited random walk.

We first recount the qualitative argument for the transition between
recurrence and transience of the $k$-excited walk that was given in
Ref.~\cite{Z04}.  Whenever the walk is in a region that contains cookies,
there is ostensibly a bias velocity $v=p-q$.  Thus the time that the walk
takes to travel a distance $L$ should be $L/(p-q)$.  If this time is less
than the total number of cookies $kL$ in this region, then the walk does not
have sufficient time to consume all the cookies as it moves.  Consequently,
the bias will persist forever and the walk will be transient.  The condition
for transience is therefore ${L}/{v}<kL$, or
\begin{eqnarray}
\label{pc}
p>p_c=\frac{1}{2}\left(1+\frac{1}{k}\right).
\end{eqnarray}
As we shall see, the behavior of the 2-excited walk is qualitatively
different in the regimes $p<p_c$ and $p>p_c$.

\subsection{Probability distributions}

\begin{figure}[ht] 
 \vspace*{0.cm}
 \includegraphics*[width=0.42\textwidth]{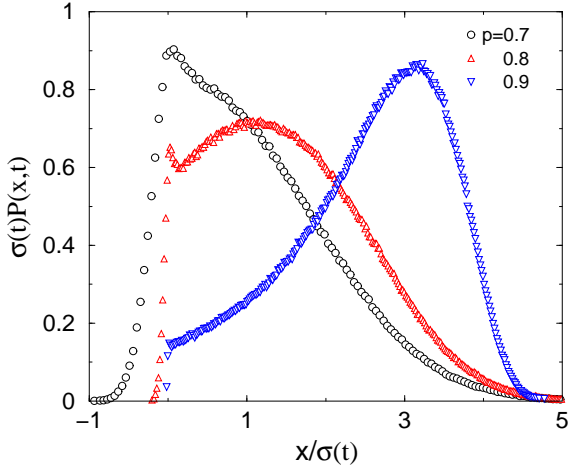}
 \caption{Scaled probability distributions of the 2-excited random walk for
   85,223 time steps and for $2\times 10^6$ realizations.  The data have been
   averaged over a $\frac{1}{2}$\% range. }
\label{pdf-2e}
\end{figure}

\begin{figure}[ht] 
 \vspace*{0.cm}
 \includegraphics*[width=0.42\textwidth]{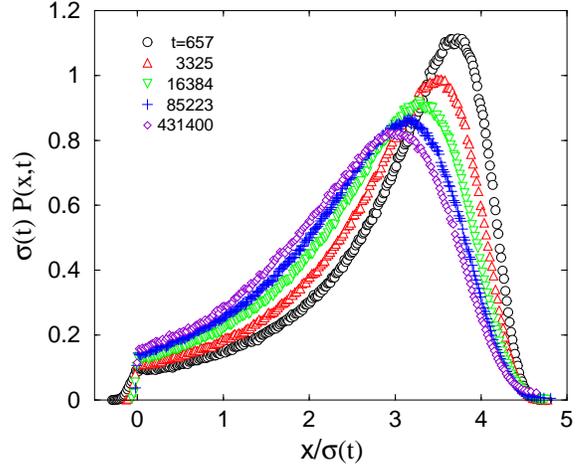}
 \caption{Scaled probability distributions of the 2-excited random walk for
   $p=0.9$ at different times for $2\times 10^6$ realizations.  The data have
   been averaged over a $\frac{1}{2}$\% range. }
\label{pdf-2c-non}
\end{figure}

Fig.~\ref{pdf-2e} shows scaled displacement probability distributions for the
2-excited random walk for representative values of $p$.  As $p$ passes
through $p_c$, there is a sudden drop in the integrated probability that the
walk lies in the range $x<0$.  Another basic feature of the distribution is
the departure from scaling, as shown in Fig.~\ref{pdf-2c-non}.  For $p=0.9$,
data collapse has not yet set in by $431,400\approx 1.5^{32}$ time steps.  If
scaling were to hold, it would only arise after an extraordinarily long time.
In the absence of such a clearly identifiable long time scale, it seems that
the asymptotic distribution will not obey single-parameter scaling.

\begin{figure}[ht] 
\vspace*{0.cm}
 \includegraphics*[width=0.4\textwidth]{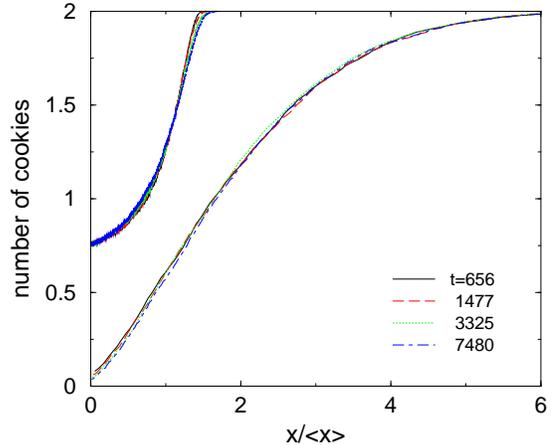}
 \caption{Cookie profile for the 2-excited random walk at $p=0.6$ (lower
   data) and at $p=0.9$ (upper data) for times 656, 1477, 3325, and 7480.
   The $x$-axis has been scaled by the average displacement of the walk at
   time $t$.}
\label{profile}
\end{figure}

The transition in the behavior of the probability distribution is mirrored in
the spatial profile of the uneaten cookies.  For $p<p_c$, the 2-excited
random walk is recurrent and thus the region close to the origin should
ultimate be devoid of cookies.  This behavior is illustrated in
Fig.~\ref{profile}, where the scaled density of uneaten cookies is shown.
For $p=0.6$ this density profile follows single-parameter scaling and the
also decays to zero near the origin because of the recurrence of the walk for
$p<3/4$.  On the other hand, for $p=0.9$, the cookie density shows a very
small departure from scaling and the density goes to a non-zero limit near
the origin as a result of the transient nature of the walk.  We can
understand this limiting cookie density near the origin by the same argument
that led to the criterion for $p_c$ in Eq.~(\ref{pc}).  For general $p>p_c$,
a walk typically resides for a time $1/(p-q)$ at each site and should eat
$1/(p-q)$ cookies, thus leaving a residue of $2-1/(p-q)$ uneaten cookies per
site.  For the case $p=0.9$ shown in Fig.~\ref{profile}, we therefore expect
3/4 of a cookie per site near the origin.

\begin{figure}[ht] 
  \vspace*{0.cm} \includegraphics*[width=0.42\textwidth]{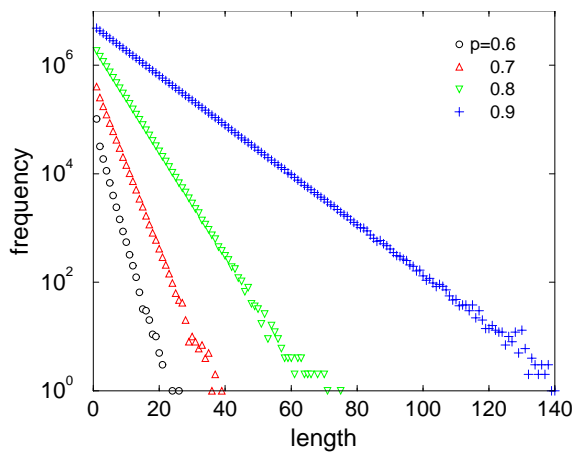}
  \caption{Length distribution of one-cookie domains on a semi-logarithmic
    scale for 50000 realizations at time $t= 37,877\approx 1.5^26$ for
    $p=0.6, 0.7, 0.8$ and 0.9.}
\label{1-domain}
\end{figure}

A more detailed characterization of the cookie profile is provided by the
size distributions of domains that contain one cookie and contain zero
cookies.  The distribution of one-cookie domains is an exponentially decaying
function of domain length (Fig.~\ref{1-domain}).  This behavior is simple to
understand.  The process that gives rise to a one-cookie domain as $p\to 1$
is a walk that makes a single step to the left after a run of consecutive
steps to the right.  During this run of rightward steps into
previously-unvisited territory (each occurring with probability $p$), one
cookie is left behind at each site.  A step to the left (which occurs with
probability $q$) terminates this one-cookie domain by a cookie-free gap.
Thus the probability of a one-cookie domain of length $k$ is simply $p^kq =
q\,e^{-k\ln(1/p)}$.  From linear least-squares fits to the data in
Fig.~\ref{1-domain}, the slopes for $p=0.6, 0.7, 0.8$, and 0.9 are,
respectively, $-0.472$, $-0.348$, $-0.219$ and $-0.104$.  On the other hand
the corresponding slopes from by the distribution $q\,e^{-k\ln(1/p)}$ are
$-0.511$, $-0.357$, $-0.223$ and $-0.105$, in excellent agreement with the
data, especially as $p\to 1$.

\begin{figure}[ht] 
 \vspace*{0.cm}
 \includegraphics*[width=0.42\textwidth]{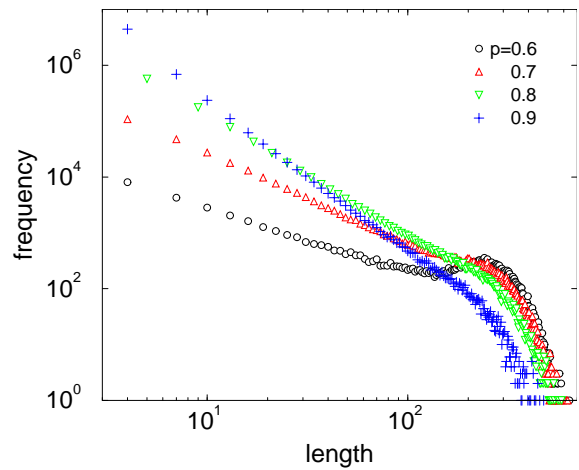}
 \caption{Length distribution of zero-cookie domains on a double logarithmic
   scale for the same parameters as Fig.~\ref{1-domain}.}
\label{0-domain}
\end{figure}

The length distribution of zero-cookie domains is more intriguing.  Such
domains arise both by the growth of isolated zero-cookie domain and by the
merging of two neighboring zero-cookie domains when an intervening one-cookie
domain disappears.  The zero-cookie domains also have qualitatively different
behavior for $p<p_c$ and $p>p_c$.  In the former case, the distribution has a
secondary peak close to the maximum-length domains, while the distribution
appears to have a pure power-law form for $p>p_c$. 

\subsection{Transit time}

Our numerical simulations indicate that the mean displacement $\langle
x(t)\rangle$ of the 2-excited walk grows with time as $t^\nu$, with $\nu$
apparently dependent on $p$.  For the 3 cases of $p=0.7$, 0.8 and 0.9,
least-squares fits to the data in the range $t>20000$ give $\nu=0.54, 0.64$,
and 0.85, respectively.  Naively, one might anticipate that the transience of
the walk for $p>3/4$ would lead to an average displacement that grows
linearly with time.  Indeed, the argument that gave the threshold bias
$p_c=3/4$ was based on the walk possessing a non-zero bias velocity
proportional to $p-q$ for $p>p_c$.  On the other hand, it has been proven
that the mean velocity of the 2-excited walk, defined as $\lim_{t\to\infty}
\langle x(t)\rangle/t$ equals zero for all $p<1$ \cite{Z04}.

While we do not know how to determine the time dependence of the
displacement, we can address the complementary problem of estimating how long
it takes a 2-excited walk to traverse a finite interval.  We already know,
from the 1-excited walk, that the mean time to traverse any cookie-free
interval within an infinite system is divergent.  Thus we now consider the
time for a 2-excited walk to traverse a finite cookie-free interval.  As we
shall show, the transit time grows faster than linearly in the interval
length when $p$ is close to but strictly less than 1.

\begin{figure}[ht] 
 \vspace*{0.cm}
 \includegraphics*[width=0.45\textwidth]{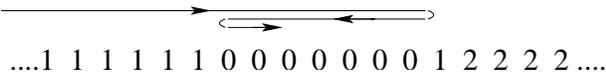}
 \caption{A 2-excited random walk in which the walk has made $k$ consecutive
   steps to the left after a long run of rightward steps.  The corresponding
   cookie configuration is also shown.}
\label{2c-disp}
\end{figure}

The 2-excited random walk starts at $x=0$ and we compute the time to reach
$x=L$, where $L\gg (1-p)^{-k}$.  This inequality ensures that a ``defect'',
namely, a region where the walk takes a small number $k$ of consecutive steps
to the left after a long string of rightward steps, will certainly occur in
the interval.  This defect represents the lowest-order perturbation to the
purely ballistic motion for the 2-excited walk when $p\to 1$.  After the
defect has been completed and a step to the right has eventually been made,
the walk finds itself at the point $x=2$ in a cookie-free region $[0,k]$.  We
ask: what is the mean time for the walk to reach the right edge of this
interval?  As illustrated in Fig.~\ref{2c-disp}, this system is a 1-excited
random walk in disguise.

We denote by $T_k$ the mean time for the 2-excited walk to reach the point
$x=k$ in this geometry.  To simplify the calculation, we assume that whenever
the left edge of the domain is reached and the last cookie is eaten, the next
step of the walk is always to the right.  Because this next step could be to
the left, our approach represents a lower bound to the transit time that we
denote by $T_k^{(\rm LB)}$.

As defined in Eq.~(\ref{exitprob_abs}) and (\ref{exittime_abs}), with
probability $E_+(k,2)$, the walk reaches $x=k$ without touching $x=0$ and the
time for this event is $t_+(k,2)$.  On the other hand, with probability
$E_-(k,2)$ the walk hits the left edge and then, according to our
approximation, hops one step to the right.  Since this second argument always
equals 2 in all formulae involving $E_\pm$ and $t_\pm$, we henceforth drop
this argument for simplicity.  After these two events, the walk is a
distance 2 from the left edge of a cookie-free domain of length $k+1$ and the
elapsed time is $t_-(k)+1$.  Thus $T_k^{(\rm LB)}$ obeys the recursion
\begin{equation}
\label{T-rec}
T_k^{(\rm LB)}=E_+(k)t_+(k)+E_-(k)\big[t_-(k)+1+T_{k+1}^{(\rm LB)}\big].
\end{equation}
If we allowed the walk to hop to the left (with probability $q$) when the
leftmost cookie is reached, then Eq.~(\ref{T-rec}) would be replaced by the
exact equation
\begin{eqnarray*}
T_k^{(\rm ex)}&=&E_+(k)t_+(k)+E_-(k)\big[t_-(k)+p(1+T_{k+1}^{(\rm ex)})\big.\nonumber\\
&+&\big.qp(2+T_{k+2}^{(\rm ex)})+q^2p(3+T_{k+3}^{(\rm ex)})+\ldots\big].
\end{eqnarray*}
The additional terms represent the contribution of paths in which the random
walk makes 1, 2, 3, {\it etc.} consecutive steps to the left upon hitting the
last cookie before a step to the right is made.  This exact equation is not
recursive and for simplicity we therefore consider the approximation
recursion for the lower bound in what follows.

By iterating (\ref{T-rec}), we express $T_k^{(\rm LB)}$ compactly as
\begin{eqnarray}
\label{T-rec1}
T_k^{(\rm LB)}=\tau_k&+&\sum_{\ell=1}^\infty \tau_{k+\ell}\prod_{j=1}^\ell E_-(k+j-1)\nonumber\\
&+&\sum_{\ell=1}^\infty \prod_{j=1}^\ell E_-(k+j-1).
\end{eqnarray}
where $\tau_k=E_-(k)t_-(k)+E_+(k)t_+(k)=2(k-2)$ is the unconditional exit
time, namely, the time for the random walk to exit either side of the
interval $[0,k]$ when the starting point is $x=2$.  Using this result for
$\tau_k$, as well as that for $E_-(k)$ in Eq.~(\ref{T-rec1}), we obtain after
straightforward algebra:
\begin{eqnarray}
\label{T-soln}
T_k^{(\rm LB)}=3(k-2)+2(k-1)(k-2)\sum_{j=k}^L \frac{1}{j}.
\end{eqnarray}
The upper limit in the sum accounts for case where a cookie-free defect of
length $k$ occurs close to the right edge of the domain $[0,L]$ so that terms
or order $L$ are needed to account for the growth of the cookie-free defect.
Thus we infer that mean time for a 2-excited walk to traverse a {\it finite}
system of length $L$ scales at least as fast as $L\ln L$.  This result is
then consistent with $\langle x(t)\rangle$ growing slower than linearly with
$t$.

\section{Summary}

The excited random walk is a very rich stochastic process in which the random
walk motion is influenced by a transitory bias.  Each site of a lattice is
initially populated with a finite number of cookies.  When a walk reaches a
site that contains cookies, the walk eats the cookie and there is a $p/q$
right/left bias at the next step.  On the other hand, the walk undergoes
unbiased motion if the visited site does not contain cookies.  Thus as the
cookies are eaten, the bias gradually disappears.  If a walk is recurrent,
then the interplay between the recurrence and the bias is particular delicate
and the excited random walk will have very different behavior from pure
diffusion or from diffusion with a constant bias.  For this reason, the
excited random walk appears to be most interesting in one dimension.

In the simplest case of the 1-excited walk on the positive half-line, we
employed a physical approach, complementary to that of Dickman et al.\ 
\cite{DA01,DAA03}, to show that the first-passage probability for the walk to
hit the origin at time $t$ decays as $t^{-1-q}$ for any $q>0$.  Thus the
1-excited walk is always recurrent, except for the trivial case of perfect
bias $p=1$.  There is especially paradoxical behavior in the case of $p<1/2$,
where the cookie-induced bias is directed toward the origin.  Here the walk
exhibits two seemingly incompatible features: the mean time for the walk to
hit the origin is infinite, while the walk strays only a finite distance from
the origin during this infinite excursion.

On the infinite line, the 1-excited walk has the same first-passage
characteristics as on the semi-infinite line.  However, the probability
distribution of the walk and the probability distributions for the positions
of the last cookie exhibit unusual features.  For the probability
distribution of the walk, there appears to be a singularity at the origin
that does not obey scaling while the rest of the distribution does scale.
There is also an asymmetry in the span of the walk that we were able to
characterize by a probabilistic argument.  The probability distributions for
the positions of the leftmost and rightmost cookies have power law tails that
we have been unable, thus far, to explain.  These tails seem to be connected
with the power-law behavior for the first-passage probability itself.

The 2-excited walk shows a wide range of strange phenomena.  A transition
between (conventional) recurrence and transience occurs as $p$ passes through
a threshold value $p_c=3/4$.  For $p<p_c$, the 2-excited walk appears to
qualitatively resemble the 1-excited walk.  On the other hand, for $p>p_c$,
the 2-excited walk is transient and this attribute appears responsible for
several unexpected features, such as the lack of scaling in the probability
distribution of the walk and in the spatial distribution of cookies, as well
as the existence of very different distributions of 1-cookie and 0-cookie
domains.  Most of these features do not yet have clean explanations and the
model is ripe for further exploration.

\acknowledgments{We thank E. Ben-Naim and P. Krapivsky for useful
  discussions.  TA gratefully acknowledges financial support from the Swiss
  National Science Foundation under the fellowship 8220-067591.  SR
  acknowledges financial support from NSF grants DMR0227670 (at BU) and DOE
  grant W-7405-ENG-36 (at LANL).  }

\appendix*
\section{Master equation}

In this appendix, we outline a Master equation approach to determine the
survival probability of the 1-excited random walk on the infinite half-line.
We merely follow the method introduced by Dickman et al.\ \cite{DA01,DAA03}
to describe an unbiased walk between a fixed absorbing boundary at $x=0$ and
a moving reflecting boundary that moves to the right with a specified
probability each time it is hit by the walk.  This perturbed random walk
motion is essentially the same mechanism as that of the excited random walk.
While the excited walk is ostensibly not Markovian, the process can be made
Markovian by considering an extended state space that is defined by the
position of the walk $x$ and also the positions of the leftmost and rightmost
cookies, $\overline y$ and $y$.  The exact master equation for the joint
probability distribution of the position of the walk and the positions of the
last cookies $P(x,\overline y, y, t)$ can be obtained in terms of a Markov
process.

We consider the simpler problem of the 1-excited random walk with an
absorbing boundary at the origin.  The system is now described by the joint
probability for the cookie to be at $x$ and the leftmost cookie to be at $y$,
$P(x,y,t)$.  This distribution obeys the recursion formula
\begin{equation}
\label{app:P}
P(x,y,t+1)=\frac{1}{2}P(x-1,y,t)+\frac{1}{2}P(x+1,y,t) ~,
\end{equation}
for $x=1,2,\ldots,y-3$.
At the origin there is an absorbing boundary condition
\begin{equation}
P(0,y,t)=0 ~,
\end{equation}
while the effect of the cookies is described by slightly modified equations
for the probability when the walk is sufficiently close to the position of
the last cookie:
\begin{eqnarray}
\label{app:Pedge}
P(y\!-\!2,y,t\!+\!1) &=& \frac{1}{2}P(y\!-\!3,y,t)\\ \!&\!+\!&\!  
\frac{1}{2}P(y\!-\!1,y,t) + 
qP(y\!-\!1,y\!-\!1,t) \nonumber\\
P(y\!-\!1,y,t\!+\!1) &=& \frac{1}{2}P(y\!-\!2,y,t)  \nonumber\\
P(y,y,t\!+\!1) &=& \frac{1}{2}P(y\!-\!1,y,t) + pP(y\!-\!1,y\!-\!1,t) \nonumber ~.
\end{eqnarray}
These equations are defined over the wedge domain $x\ge0, x\le y$.  We choose
the initial conditions that the walk starts at $x=1$ and that the last cookie
is also at $x=1$; thus $P(x,y,0)=\delta_{x,1}\delta_{y,1}$.  These equations
can be iterated then numerically and this is the method that we used for the
simulations presented in Sec.~II.

{}From these Master equations, we compute the probability for the excited
random walk to reach the origin for the first time at time $t$, from which
the survival probability of the walk follows.  We merely outline this
derivation as it very similar to that given by Dickman et al.\ 
\cite{DA01,DAA03}.  We first rewrite the Master equations in terms of the
generating function $\hat P(x,y,z)=\sum_{t=0}^\infty z^t P(x,y,t)$.  In this
representation, the solution can be written in the product form
\begin{equation}
\hat P(x,y,z)= \hat A(x,z) \hat B(y,z),
\end{equation}
where
\begin{equation}
\hat A(x,z)=\lambda^x-\lambda^{-x}\,,\quad{\rm with}\quad
\lambda=\frac{1}{z}+\sqrt{\frac{1}{z^2}-1},
\end{equation}
and where $\hat B$ is defined recursively by
\begin{equation}
\label{hatB}
\frac{\hat B(y+1,z)}{\hat B(y,z)} = \frac{\tanh \Lambda y + 2q\Lambda}{\tanh \Lambda 
y + \Lambda} ~,
\end{equation}
in the $\epsilon\equiv z-1\to 0$ limit, with $\Lambda=\ln\lambda$.

The long-time behavior of the survival probability $\hat
S(z)=\sum_{x,y}P(x,y,z)$ can now be obtained from the $z\to 1$ limit
($\Lambda\to 0$) of $P(x,y,z)$.  For any finite $y$, and thus for $y=y_0\sim
2(p-1)$, $\hat B(y_0) \sim \hat B(1) \sim 1/\Lambda$ as $\Lambda\to 0$.  On
the other hand, for $y>y_0$ we approximate the recursion formula (\ref{hatB})
by
\begin{eqnarray}
 \ln \frac{\hat B(y)}{\hat B(y_0)} &=& \sum_{k=y_0}^{y-1} \ln 
 \frac{1+2(p-1)\Lambda/\tanh\Lambda k}{1+\Lambda/\tanh\Lambda k}\nonumber\\
&\approx& (2p-3)\Lambda\sum_{k=y_0}^{y-1}  \frac{1}{\tanh{\Lambda k}}\\
&\approx& (3-2p) \ln \frac{\sinh \Lambda y_0}{\sinh \Lambda y}\nonumber
\end{eqnarray}
Putting these terms together we obtain
\begin{eqnarray}
\hat S(z) &=& \sum_{y=1}^\infty \hat B(y,z) \sum_{x=0}^y \hat A(x,z)\nonumber\\
&\sim& \sum_{y=1}^\infty 
\hat B(y,z)\frac{4}{\Lambda}\sinh^2\frac{\Lambda y}{2}\nonumber\\
&\sim& \Lambda^{1-2p}\sum_{y=y_0}^\infty \frac{\sinh^2\Lambda 
y/2}{\sinh^{3-2p}\Lambda y}\\
&\sim& \Lambda^{-2p} \int_{\Lambda y_0}^\infty \frac{\sinh^2 
u/2}{\sinh^{3-2p} u} du ~. \nonumber
\end{eqnarray}
Since the integral above goes to a finite constant as $\Lambda\to 0$, $\hat
S(z) \sim \Lambda^{-2p} \sim \epsilon^{-p}$, which corresponds to $S(t)\sim
t^{-q}$.

For the 1-excited random walk on the infinite line the phase space is now the
three-dimensional domain defined by $\overline y \le x\le y$.  In principle,
we can write recursion formulae similar to Eqs.~(\ref{app:P}) --
(\ref{app:Pedge}) for $P(x,\overline y, y, t)$, but the solution to these
equations no longer seems practical.


\begin{thebibliography}{99}
  
\bibitem{PW97} M. Perman and W. Werner, Probab.\ Theory Related Fields {\bf
    108}, 357-383 (1997).

\bibitem{D99} B. Davis, Probab.\ Theory Related Fields {\bf 113}, 501-518 (1999).
  
\bibitem{P01} R. Pemantle, {\it Random processes with reinforcement},
  preprint.  Available from:

 www.math.upenn.edu/$\sim$pemantle/papers/Papers.html.

\bibitem{BW03} I. Benjamini and D. B. Wilson, Elect.\ Comm.\ in Probab. {\bf
    8}, 86-92 (2003).

\bibitem{ABV03} O. Angel, I. Benjamini and B. Vor\'ag, Elect.\ Comm.\ in Probab. {\bf
    8}, 6-16 (2003).

\bibitem{Z04} M. P. W. Zerner, {\it math.PR/0403060}. 
  
\bibitem{SAW} N. Madras and G. Slade, {\it The Self-Avoiding Walk}
  (Birkhauser, Boston, 1992).

\bibitem{TSAW} D. J. Amit, G. Parisi, and L. Peliti, Phys.\ Rev.\ B {\bf 27},
  1635 (1983); B. T\'oth, J. of Stat.\ Phys.\, {\bf 77}, (1994); {\it ibid.,\ } 
The Annals of Probability {\bf 23}, 1523 (1995).

\bibitem{PV99} R. Pemantle, Probab.\ Th.\ and Rel.\ Fields, {\bf 92}, 117
  (1992); R. Pemantle and S. Volkov, Ann.\ Probab.\ {\bf 27}, 1368 (1999).

\bibitem{DA01} R. Dickman and D. ben-Avraham, Phys.\ Rev.\ E {\bf 64},
  020102(R) (2001).

\bibitem{DAA03} R. Dickman, F. F. Araujo, Jr., and D. ben-Avraham, Brazilian
  J. Phys.\ {\bf 33}, 450 (2003).

\bibitem{F68} W. Feller {\it An Introduction to Probability Theory and Its
Applications,} (Wiley, New York, 1968).

\bibitem{fpp} S. Redner, {\it A Guide to First-Passage Processes}, (Cambridge
University Press, New York, 2001).

\bibitem{AHZ} Similar arguments were used in T. Antal, H. J. Hilhorst, and R.
  K. P. Zia, J. Phys.\ A {\bf 35}, 8145 (2002).

\bibitem{concmat} R. L. Graham, D. E. Knuth, and O. Patashnik, {\it Concrete
    Mathematics: A Foundation for Computer Science} (Addison-Wesley, Reading,
  MA, 1989).

\bibitem{AS} M. Abramowitz and I. A. Stegun, {\it Handbook of Mathematical
    Functions} (Dover Publications Inc., New York, 1970).

\bibitem{knuth} R.~L.~Graham, D.~E.~Knuth, and O.~Patashnik, {\em Concrete
    Mathematics: A Foundation for Computer Science} (Reading, Mass.:
  Addison-Wesley, 1989).

\bibitem{W} G. H. Weiss, {\it Aspects and Applications of the Random Walk},
(North-Holland, Amsterdam, 1994).

\bibitem{GR} I. S. Gradshteyn and I. M Ryzhik, {\it Tables of Integrals,
    Series, and Products}, (Academic Press, New York, 1965).

\end{thebibliography}
\end{document}